\documentclass[reqno]{amsart}
\usepackage{fullpage}
\usepackage{bm}
\usepackage{amsfonts,mathrsfs,graphicx,color}
\usepackage[latin1]{inputenc}
\usepackage[english]{babel}

\newtheorem{theorem}{Theorem}

\newtheorem{lemma}{Lemma}

\newtheorem{conjecture}{Conjecture}

\title{ \bf Polytopal complexes: maps, chain complexes and... necklaces}
\author{Frédéric Meunier}
\address{Université Paris Est, LVMT, Ecole Nationale des Ponts et Chaussées, 6-8 avenue
Blaise Pascal, Cité Descartes Champs-sur-Marne, 77455
Marne-la-Vallée cedex 2, France.} \email{ frederic.meunier@enpc.fr}
\date{\today}

\begin{document}

\begin{abstract}
The notion of {\em polytopal map} between two polytopal complexes is
defined. Surprisingly, this definition is quite simple and extends
naturally those of simplicial and cubical maps. It is then possible
to define an induced chain map between the associated chain
complexes. Finally, we use this new tool to give the first
combinatorial proof of the splitting necklace theorem of Alon. The
paper ends with open questions, such as the existence of Sperner's lemma
for a polytopal complex or the existence of a cubical approximation
theorem.
\end{abstract}

\maketitle

\section*{Introduction}

At the very beginning of algebraic topology, there is a simple
construction that allows to derive from a simplicial map a chain map
that carries a formal sum of oriented simplices onto a formal sum of
their image. This chain map can then be treated as an algebraic
object and the powerful machinery of algebra can be used to derive
properties of the chain map, and then of the starting simplicial map
itself. Ehrenborg and Heytei have shown in 1995 \cite{EhHe95} that a
similar construction is possible for cubical maps, that is the
analog of simplicial maps for cubical complexes. A natural question
is then to ask whether it is possible to define in a similar way
{\em polytopal maps} between two polytopal complexes and to define
naturally an associated chain map that would keep the essential
properties of the polytopal map. To our knowledge, this definition
is something new.

The purpose of the present work is to show that it is possible. We
give a very natural definition of a polytopal map and we will see
that this definition contains both the notions of simplicial and
cubical maps, although our definition of a cubical map is then more
restrictive than which of Ehrenborg and Hetyei. This point will be
discussed. Then we show how to derive a chain map from it. Finally,
we use this new framework to give the first purely combinatorial
proof of the celebrated splitting necklaces theorem, which states
that an open necklace with $t$ types of beads and a multiple of $q$
of each type can always be fairly divided between $q$ thieves using
no more than $t(q-1)$ cuts. According to Ziegler \cite{Zi02}, a
combinatorial proof in a topological context is a proof using no
simplicial approximation, no homology, no continuous map.

The plan of the paper is the following. In the first section
(Section \ref{sec:pol}), we recall the definition of a polytopal
complex (Subsection \ref{subsec:elementary}), we define the notion
of polytopal map (Subsection \ref{subsec:polymap}), the notion of
oriented polytope and the notion of boundary operator for a
polytopal complex, inducing the notion of chain complex (Subsection
\ref{subsec:chaincomplex}). It is possible to define the cartesian
product of two chains - a thing that is not possible for simplicial
complexes (Subsection \ref{subsec:product}). Then we will see how a
polytopal map induces a chain map. (Subsection
\ref{subsec:chainmap}). In Subsection \ref{subsec:homotopy}, an
homotopy equivalence of polytopal maps is proved. In the second
section (Section \ref{sec:necklace}), we prove combinatorially the
splitting necklace theorem. Actually, we get also a direct proof of
the generalization found by Alon, Moshkovitz and Safra
\cite{AlMoSa06} when there is not necessarily a multiple of $q$
beads of each type. In the last section (Section
\ref{sec:discussion}), open questions are presented (Sperner's lemma
for a polytopal complex, cubical approximation theorem).

\medskip

{\bf Acknowledgement}: The author thanks Robin Chapman for his valuable comments, in particular his comment on a definition of a polytopal map, given in a first version of this paper, that was not a true generalization of the notion of simplicial map.

\section{Polytopal and chain maps}

\label{sec:pol}

\subsection{Polytopal complexes}

\label{subsec:elementary}

A {\em polytopal complex} $\mathsf{P}$ is a collection of finite
sets (called {\em faces}) on a vertex set $V(\mathsf{P})$ such that
\begin{enumerate}
\item For every $\sigma\in\mathsf{P}$ the elements of $\sigma$ can
be represented as the vertices of a finite dimensional polytope,
where the faces contained in $\sigma$ are exactly the vertex sets of
the faces of this polytope,
\item If $\sigma,\tau\in\mathsf{P}$ then $\sigma\cap\tau=\emptyset$
or $\sigma\cap\tau\in\mathsf{P}$.
\end{enumerate}

For every face $\sigma$ the dimension of the polytope
associated to $\sigma$ is called the dimension of $\sigma$ and denote it by
$\dim\sigma$. Any face $\tau$ of $\sigma$ such that
$\dim\tau=\dim\sigma-1$ is called a facet of $\sigma$. When all
polytopes are simplices, the polytopal complex is called a {\em
simplicial complex}, and when all polytopes are cubes, it is called
a {\em cubical complex}. When $\mathsf{P}$ and $\mathsf{Q}$ are two
polytopal complexes, $\mathsf{P}\times\mathsf{Q}$ is the
$(\dim\mathsf{P}+\dim\mathsf{Q})$-polytopal complex whose faces are
all $\sigma\times\tau$ where $\sigma\in\mathsf{P}$ and
$\tau\in\mathsf{Q}$. When the cartesian product of
$\mathsf{P}$ is taken $s$ times by itself, we use the notation
$\mathsf{P}^s$.

Finally, we emphasize the fact that we assume that there is no
degeneracy in the sense that the faces of the polytopes are always
true faces: two distinct faces of the same polytope have distinct
supporting affine subspaces.

\subsection{Polytopal maps}

\label{subsec:polymap}

We give now the definition of a polytopal map.

Let $\mathsf{P}$ and $\mathsf{Q}$ be two polytopal complexes. A map
$\lambda: V(\mathsf{P})\rightarrow V(\mathsf{Q})$ is a {\em
polytopal map} if for every $d$-face $\sigma$ of $\mathsf{P}$, $\lambda(\sigma)$ is a subset of a $d'$-face of
$\mathsf{Q}$, with $d'\leq d$.

\begin{figure}
\includegraphics[width=10cm]{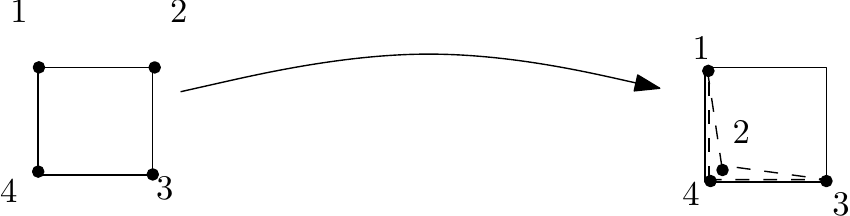}
\caption{A cubical map.} \label{fig:cubicalmap}
\end{figure}

In Figure \ref{fig:cubicalmap}, a polytopal map in the particular
case of cubical complexes is illustrated. Note that the image of a
cube is not necessarily a cube.

When $\mathsf{P}$ and $\mathsf{Q}$ are two simplicial complexes, the
classical definition of a simplicial map is the following: a map
$\lambda: V(\mathsf{P})\rightarrow V(\mathsf{Q})$ is a simplicial
map if for every simplex $\sigma$ of $\mathsf{P}$, $\lambda(\sigma)$
is a simplex of $\mathsf{Q}$. It is straightforward to check that,
in this case, this definition coincides with the one given above.

When $\mathsf{P}$ and $\mathsf{Q}$ are two cubical complexes, the
classical definition of a cubical map is the following (see Fan
\cite{Fa60} or Ehrenborg and Hetyei \cite{EhHe95} for instance): a
map $\lambda: V(\mathsf{P})\rightarrow V(\mathsf{Q})$ is a {\em
cubical map} if the following conditions are both fulfilled:
\begin{enumerate}
\item for every cube $\sigma$ of $\mathsf{P}$, $\lambda(\sigma)$ is a subset of the vertices of a cube of $\mathsf{Q}$
\item $\lambda$ takes adjacent vertices to adjacent vertices or the same vertex.
\end{enumerate}

Our polytopal map when the polytopal complexes are cubical complexes
is a cubical map in the sense of Fan or Ehrenborg and Hetyei, but
the converse is not true as the following example shows\footnote{In
a first version of this work, presented at the TGGT 08 conference in
Paris, the author was not conscious of this fact and claimed
erroneously that his notion of a polytopal map contains the notion
of a cubical map by Fan, Ehrenborg and Hetyei.}. Indeed, the
following map between a $3$-dimensional cube and a $4$-dimensional
cube is a cubical map in the sense Fan-Ehrenborg-Hetyei but not in
the sense of the present paper. The minimal face containing the
image of the $3$-cube is the $4$-cube itself.

$$\begin{array}{ccc}
(0,0,0) & \rightarrow & (0,0,0,0) \\
(1,0,0) & \rightarrow & (1,0,0,0) \\
(0,1,0) & \rightarrow & (0,1,0,0) \\
(1,1,0) & \rightarrow & (0,0,0,0) \\
(0,0,1) & \rightarrow & (0,0,0,1) \\
(1,0,1) & \rightarrow & (0,0,0,0) \\
(0,1,1) & \rightarrow & (0,0,0,0) \\
(1,1,1) & \rightarrow & (0,0,1,0) \\
\end{array}$$

Whether it is possible to define a more general notion of a
polytopal map which would contain a cubical map in this more general
sense is an open question (see Section \ref{sec:discussion} for a
complementary discussion).

\subsection{Oriented polytopes and chain complex of polytopal
complexes}

\label{subsec:chaincomplex}

Let $P$ be a $d$-dimensional polytope. $\epsilon:
{V(P)\choose d+1}\rightarrow\{-1,0,+1\}$ is an {\em orientation} of
$P$ if the following subsequent conditions are fulfilled:

\begin{enumerate}
\item $\epsilon(v_0,\ldots,v_d)=0$ if and only if $v_0,\ldots,v_d$ are
affinely dependent.
\item
$\epsilon(v_0,\ldots,v'_i,\ldots,v_d)=\epsilon(v_0,\ldots,v_i,\ldots,v_d)$
if $v_i$ and $v'_i$ are in the same open half-space (of
$\mathbb{R}^d$) delimited by the supporting hyperplane of
$v_0,\ldots,\hat v_i,\ldots,v_d$, and
$\epsilon(v_0,\ldots,v'_i,\ldots,v_d)=-\epsilon(v_0,\ldots,v_i,\ldots,v_d)$
if not (the hat means a missing element).
\item
$\epsilon(v_0,v_1,\ldots,v_d)=\mbox{sign}(\pi)\epsilon(v_{\pi(0)},v_{\pi(1)},\ldots,v_{\pi(d)})$
for every permutation $\pi$.
\end{enumerate}

\begin{lemma}
Every $d$-dimensional polytope admits exactly two distinct
orientations. Moreover, if $\epsilon$ is one of them, $-\epsilon$ is
the other one.
\end{lemma}

\begin{proof}
One can define such an $\epsilon$ by the sign of a determinant, it
is enough to prove that an $\epsilon$ is completely defined by its
value on a $d$-simplex. But this is obvious since it is possible to
go from a $d$-simplex to any other simplex of the polytope by
replacing one vertex after the other, while keeping the affinely
independence.
\end{proof}

For every $d$-face $\sigma$ of a polytopal complex $\mathsf{P}$, an
{\em orientation} $\epsilon$ of $\sigma$ is an orientation of the
associated polytope. We say then that $(\sigma,\epsilon)$ is an {\em
oriented face} of $\mathsf{P}$. When $\tau$ is a facet of $\sigma$,
one defines the {\em induced orientation} of an orientation
$\epsilon$ of $\sigma$ by
$$\epsilon|_{\tau}(v_0,\ldots,v_{d-1}):=\epsilon(v_0,\ldots,v_{d-1},v_d)$$
for any $v_1,\ldots,v_{d-1}\in\tau$ and $v_d\in\sigma\setminus\tau$.

One can now define the chain complex of a polytopal complex
$\mathsf{P}$. The group of formal sums of oriented $d$-faces with
coefficients in $\mathbb{Z}$ is denoted by $C_d(\mathsf{P})$, where
$-(\sigma,\epsilon)=(\sigma,-\epsilon)$. The boundary
operator is defined on an oriented $d$-face $(\sigma,\epsilon)$:
$$\partial(\sigma,\epsilon):=\sum_{\tau\mbox{\small{ facets of
}}\sigma}(\tau,\epsilon|_{\tau}).$$

\begin{lemma}
$\partial\circ\partial=0$.
\end{lemma}

\begin{proof}
Let $\eta$ be a $(d-2)$-face of $\sigma$, and let $\tau_1$ and
$\tau_2$ the two facets of $\sigma$ containing $\eta$. The lemma
will be a consequence of the equality
$\left(\epsilon|_{\tau_1}\right)|_{\eta}=-\left(\epsilon|_{\tau_2}\right)|_{\eta}$.

This equality is obvious: take $v_0,v_1,\ldots, v_{d-2}$ be $d-2$
affinely independent vertices of $\eta$, take moreover $w_1$ (resp.
$w_2$) be an independent vertex of $\tau_1$ (resp. $\tau_2$) and $v$
another independent vertex of $\sigma$; then
$\epsilon(v_0,\ldots,v_{d-2},w_1,v)=-\epsilon(v_0,\ldots,v_{d-2},w_2,v)$.
\end{proof}

Hence, the $C_d(\mathsf{P})$'s provide a chain complex $\mathcal{C}(\mathsf{P})$.

\subsection{Product of chains}\label{subsec:product}

An interesting property of the polytopal complexes is that the cartesian product of two polytopal
complexes is still a polytopal complex. It is possible to exploit this property on the level of chains identifying
$C_d(\mathsf{K})\otimes C_{d'}(\mathsf{L})$ and $C_{d+d'}(\mathsf{K}\times\mathsf{L})$. It will be used in Subsections \ref{subsec:homotopy} and \ref{subsec:combdold}

\medskip

\noindent{\bf Remark}: This is a combinatorial interpretation of the classical {\em tensor product} of chain complexes (see p. 338 of the book by Munkres \cite{Mu84}).

\medskip

We show how this identification works for two oriented polytopes.
The extension for all chains is done through bilinearity. Let
$(\sigma,\epsilon)$ and $(\sigma',\epsilon')$ be two oriented
polytopes, the first one in $\mathsf{K}$ and the second one in
$\mathsf{L}$. Let $v_0,v_1,\ldots,v_d$ (resp.
$v'_0,v'_1,\ldots,v'_{d'}$) be $d+1$ independent vertices of $\sigma$
(resp. $d'+1$ independent vertices of $\sigma'$). Note that all pairs $(v_i,v'_j)$ are vertices of
the polytope $\sigma\times\sigma'$. Now let
$\epsilon\times\epsilon'$ be the orientation of
$\sigma\times\sigma'$ such that
$$(\epsilon\times\epsilon')((v_0,v'_0),(v_1,v'_0),\ldots,(v_d,v'_0),(v_d,v'_1),\ldots,(v_d,v'_{d'})):=\epsilon(v_0,v_1,\ldots,v_d)\epsilon'(v'_0,v'_1,\ldots,v'_d).$$

\begin{lemma}
$\epsilon\times\epsilon'$ is a well defined orientation of $\sigma\times\sigma'$,
in the sense that it is independent of the choice of the independent vertices of $\sigma$ and $\sigma'$.
\end{lemma}

\begin{proof}[Sketch of proof] A simple way to see it is through the determinant
representation. A simple computation leads to the conclusion.
\end{proof}

Define
$$(\sigma,\epsilon)\otimes(\sigma',\epsilon'):=(\sigma\times\sigma',\epsilon\times\epsilon').$$

A useful relation is given by the following lemma.

\begin{lemma}\label{lem:diffprod}
Let $c\in C_d(\mathsf{K})$ and $c'\in C_{d'}(\mathsf{L})$.
Then $$\partial c\otimes c'=(-1)^{d'}\partial c\otimes c'+c\otimes\partial c'.$$
\end{lemma}

\begin{proof}[Sketch of proof] Let $(\sigma,\epsilon)$ (resp. $(\sigma',\epsilon')$)
be an oriented $d$-dimensional polytope of $\mathsf{K}$ (resp. an oriented $d'$-dimensional polytope $\mathsf{L}$).
Take a facet $\tau$ of $\sigma$ and a facet $\tau'$ of $\sigma'$.
Compute the orientation induced on $\tau\times\sigma'$ and on $\sigma\times\tau'$ by $\epsilon\times\epsilon'$ to conclude.
\end{proof}

\subsection{Chain maps}

\label{subsec:chainmap} A polytopal map induces naturally a chain
map of the corresponding chain complexes. Recall that a chain map is
such that $\lambda_{\#}\circ\partial=\partial\circ\lambda_{\#}$.

\begin{theorem}\label{thm:chainmap}
Let $\lambda: V(\mathsf{P})\rightarrow V(\mathsf{Q})$ be a polytopal
map. Then it is possible to build a chain map $\lambda_{\#}:
\mathcal{C}(\mathsf{P})\rightarrow\mathcal{C}(\mathsf{Q})$ such that
\begin{itemize}
\item $\lambda_{\#}(v,\epsilon)=(\lambda(v),\epsilon')$ for every $v\in V(\mathsf{P})$,
with $\epsilon'(\lambda(v))=\epsilon(v)$ and
\item for every oriented $d$-face $(\sigma,\epsilon)$ of $\mathsf{P}$
and any oriented $d'$-face $(\sigma',\epsilon')$ of $\mathsf{Q}$
such that $\lambda(\sigma)\subseteq\sigma'$ and $d'\leq d$, either $d'<d$ and $\lambda_{\#}(\sigma,\epsilon)=0$, or $d'=d$ and there exists an $\alpha_{\sigma}\in\mathbb{N}$
such that
$\lambda_{\#}(\sigma,\epsilon)=\alpha_{\sigma}(\sigma',\epsilon')$.
\end{itemize}
Moreover, this construction is functorial: $\mbox{\textup{id}}$
induces identity on the level of chain complexes, and $(\lambda\circ
\mu)_{\#}=\lambda_{\#}\circ \mu_{\#}$.
\end{theorem}

\begin{proof}Let us postpone the discussion of the functoriality.
The proof works by induction on $d$.

For $d=0$, there is nothing to prove.

Suppose $d\geq 1$ and let $(\sigma,\epsilon)$ be an oriented
$d$-face of $\mathsf{P}$. Let $\sigma'$ be a $d'$-face of
$\mathsf{Q}$ such that $\lambda(\sigma)\subseteq\sigma'$ and $d'\leq d$.

A facet $\tau$ of $\sigma$ is {\em vanishing}
if $\lambda(\tau)$ is included in a face of $\sigma'$ of dimension at most $d-2$\footnote{By convention, a $-1$-face is the empty set,
in the case when $d=1$.}. 
By induction, if $\tau$ is vanishing, then $\alpha_{\tau}=0$. Note that if $\tau$ is non-vanishing, then there is a unique facet of $\sigma'$ containing $\lambda(\tau)$. Denote by $f(\tau)$ this facet.

By induction, $\lambda_{\#}$ is defined for
$\partial(\sigma,\epsilon)$: one has
$$\lambda_{\#}\left(\partial(\sigma,\epsilon)\right)=\sum_{\tau\mbox{\small{ non-vanishing facet of
}}\sigma}\alpha_{\tau}(f(\tau),\epsilon(\tau)),$$ where
$\epsilon(\tau)$ is a certain orientation of $f(\tau)$. Choose
the $\epsilon(\tau)$ such that if $f(\tau_1)=f(\tau_2)$
then $\epsilon(\tau_1)=\epsilon(\tau_2)$. Then, regrouping the
$\tau$ having same image, one defines the
$\beta_{\tau'}\in\mathbb{N}$ associated to the facets $\tau'$ of
$\sigma'$, such that
\begin{equation}\label{eq:kkk}
\lambda_{\#}\left(\partial(\sigma,\epsilon)\right)=\sum_{\tau'\mbox{\small{
facet of
}}\sigma'}\beta_{\tau'}(\tau',\epsilon(\tau')),\end{equation} where
$\epsilon(\tau')$ is a certain orientation of $\tau'$.

By induction, one has also
$$(\partial\circ\lambda_{\#})\left(\partial(\sigma,\epsilon)\right)=0$$
since $\partial$ and $\lambda_{\#}$ commute when applied to a
$(d-1)$-chain.

Hence, $$\partial\sum_{\tau'\mbox{\small{ facet of
}}\sigma'}\beta_{\tau'}(\tau',\epsilon(\tau'))=0, $$ which implies
that, for any two facets $\tau'_1$ and $\tau'_2$ sharing a common
$(d-2)$-face, one has $\beta_{\tau'_1}=\beta_{\tau'_2}$
and the induced orientations of $\epsilon(\tau'_1)$ and
$\epsilon(\tau'_2)$ on the common $(d-2)$-face are opposed.
Therefore, all the $\beta_{\tau'}$ are equal\footnote{As noted by Robin Chapman, the equality of the $\beta_{\tau'}$ is also a direct consequence of acyclicity, but we want to keep the combinatorial track.}  -- call this common
value $\alpha_{\sigma}$ -- and the $\epsilon(\tau')$ are induced by
a common orientation $\epsilon'$ of $\sigma'$. Then define
$\lambda_{\#}(\sigma,\epsilon)$ to be
$\alpha_{\sigma}(\sigma',\epsilon')$.

It remains to check that one has
$(\lambda_{\#}\circ\partial)(\sigma,\epsilon)=(\partial\circ\lambda_{\#})(\sigma,\epsilon)$.
By definition of $\epsilon'$ and $\alpha_{\sigma}$ one rewrites
Equation (\ref{eq:kkk}):
$$(\lambda_{\#}\circ\partial)(\sigma,\epsilon)=\sum_{\tau'\mbox{\small{
facet of }}\sigma'}\alpha_{\sigma}(\tau',\epsilon'|_{\tau'})).$$

Since $\sum_{\tau'\mbox{\small{ facet of
}}\sigma'}(\tau',\epsilon'|_{\tau'}))$ is precisely
$\partial(\sigma',\epsilon')$, one has
$$(\lambda_{\#}\circ\partial)(\sigma,\epsilon)=\partial\left(\lambda_{\#}(\sigma,\epsilon)\right).$$

Now, it remains to show that the identity as a polytopal map induces
a identity as chain map, but it is straightforward by following the
induction scheme above. The same holds for the composition.
\end{proof}

\begin{figure}
\includegraphics[width=13cm]{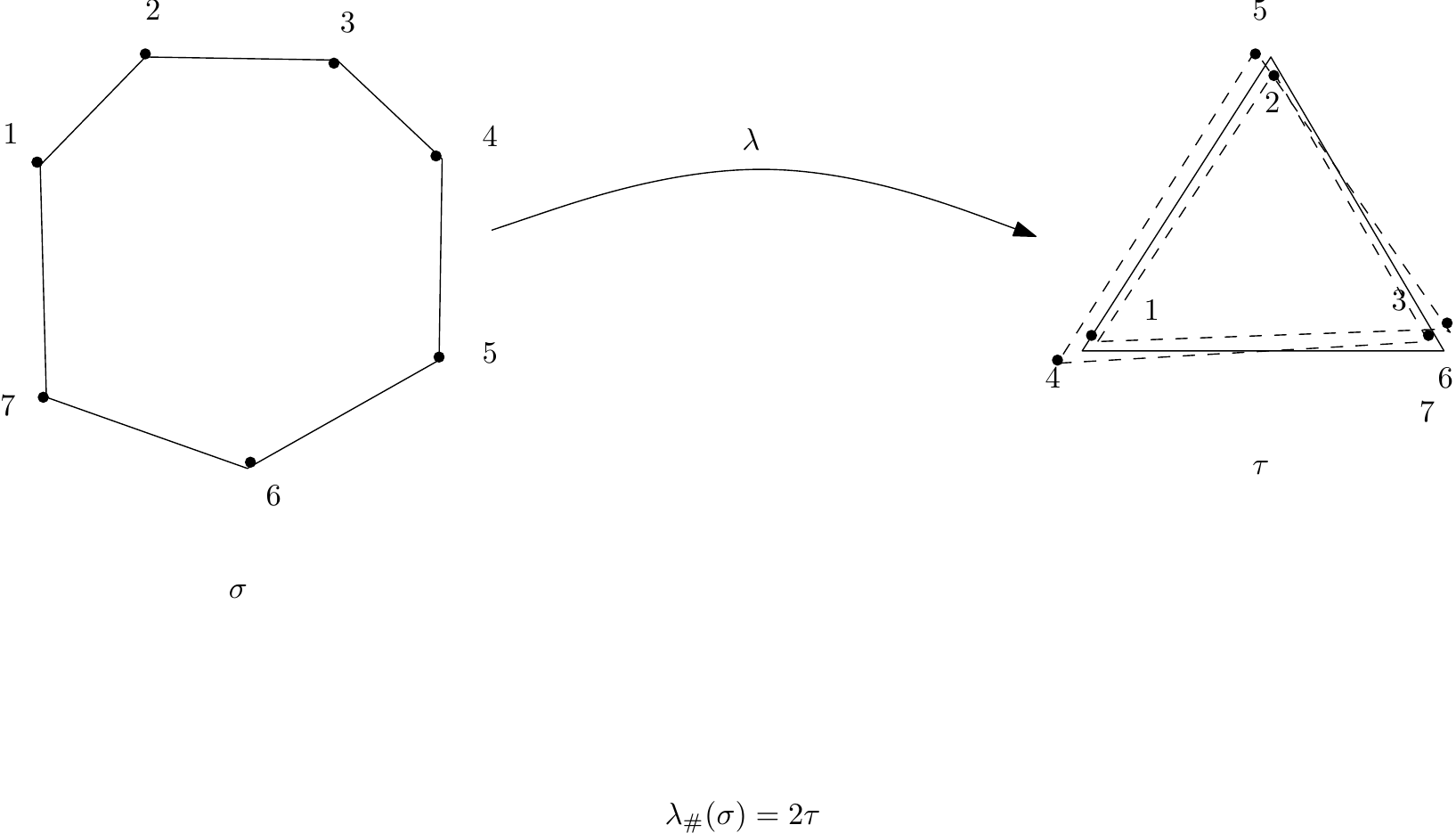}
\caption{A polytopal map with $\alpha_{\sigma}=2$}
\label{fig:polytopalmap}
\end{figure}

In Figure \ref{fig:polytopalmap}, an example of a polytopal map with
$\alpha_{\sigma}=2$ is given.

\medskip

\noindent{\bf Remark}: Note the similarity with the notion of cellular map between
CW-complexes: a polytopal map maps a polytope to a polytope of
smaller dimension, and the $\alpha_{\sigma}$ is reminiscent of the
degree of a cellular map. Nevertheless, a polytopal map is something
purely combinatorial, and the image of a polytope is not necessarily
a polytope.

\medskip

When one works with simplicial maps, the $\lambda_{\#}$ is
necessarily the classical one. Indeed, one can easily prove by
induction that the $\alpha_{\sigma}$ above is equal to $-1$, 0 or
$+1$. The same holds for cubical maps. Moreover, given a polytopal
map $\lambda$, the definition of a chain map and Theorem
\ref{thm:chainmap} above are enough to define inductively
$\lambda_{\#}(\sigma,\epsilon)$. Start by the vertices, then
define $\lambda_{\#}$ on the edges, then on the 2-faces, and so
on... 

Finally, one can show that if $\lambda(\sigma)\neq\sigma'$
(with the notation used in the theorem above), then
$\lambda_{\#}(\sigma,\epsilon)=0$ (but the converse is not
necessarily true).

\subsection{Homotopy equivalence}\label{subsec:homotopy}

An important notion when dealing with induced homology maps is the
notion of {\em homotopic maps}. Two polytopal maps
$\lambda,\mu:\mathsf{K}\rightarrow\mathsf{L}$ are {\em homotopic}
when there is a path $P_n=v_0\ldots v_n$ (seen as a $1$-dimensional
cubical map) and a polytopal map $\phi:\mathsf{K}\times
P_n\rightarrow\mathsf{L}$ such that for every $v\in V(\mathsf{K})$
we have $\lambda(v)=\phi(v,v_0)$ and $\mu(v)=\phi(v,v_n)$. If we can
take $P_n$ to be a path of length one, i.e. a standard $1$-cube
$\mathsf{\square}_1$, then we call $\lambda$ and $\mu$ {\em elementary
homotopic} maps.

Take Definition 30, page 285, in the paper by Ehrenborg and Heytei
\cite{EhHe95} to see that this notion contains the notion of
homotopic cubical maps. It contains also the notion of homotopic
simplicial maps. In this context, two ``elementary homotopic maps''
$\lambda$ and $\mu$ are so-called {\em contiguous} simplicial maps:
for each simplex $v_0,\ldots,v_d$ of $\mathsf{K}$, the points
$$\lambda(v_0),\ldots,\lambda(v_d),\mu(v_0),\ldots,\mu(v_d)$$ span a
simplex of $\mathsf{L}$. Starting with $\lambda$ and substituting
progressively the image by $\lambda$ of the vertices of $\mathsf{K}$
by their image by $\mu$, we see that two contiguous simplicial
maps are homotopic polytopal maps. Actually, the length of this path
is at most the chromatic number of the $1$-skeleton of $\mathsf{K}$
since, at each step, we can substitute the image of a stable of this
graph.

We prove now that homotopic polytopal maps induce homotopic chain
maps, that is, there is a morphism $D:C_i(\mathsf{K})\rightarrow
C_{i+1}(\mathsf{L})$ such that
$$D\circ\partial+\partial\circ D=\lambda_{\#}-\mu_{\#}.$$

\begin{lemma}
If the polytopal maps $\lambda,\mu:\mathsf{K}\rightarrow\mathsf{L}$
are homotopic, then the induced chain maps $\lambda_{\#},\mu_{\#}$
are chain homotopic.
\end{lemma}

\begin{proof}
By transitivity, it is enough to prove it for elementary homotopic
polytopal map. We denote by $v_0$ and $v_1$ the two vertices of $P_1$,
which we identify with the 1-dimensional oriented simplex $[v_0,v_1]$ (standard notation for oriented simplices).

Let $s:\mathcal{C}(\mathsf{K})\rightarrow\mathcal{C}(\mathsf{\mathsf{K}}\times P_1)$ defined by
$s(\sigma,\epsilon)=(\sigma,\epsilon)\otimes[v_0,v_1]$.

Using the definition of $s$, we compute with the help of Lemma \ref{lem:diffprod}
$$(s\circ\partial)(\sigma,\epsilon)=-(\partial(\sigma,\epsilon))\otimes[v_0,v_1].$$
On the other hand, one has
$$(\partial\circ s)(\sigma,\epsilon)=(\partial(\sigma,\epsilon))\otimes[v_0,v_1]+(\sigma,\epsilon)\otimes v_1-(\sigma,\epsilon)\otimes v_0.$$
Hence, $$(s\circ\partial+\partial\circ
s)(\sigma,\epsilon)=(\sigma,\epsilon)\otimes v_1-(\sigma,\epsilon)\otimes v_0.$$
Defining $D:=\phi\circ s$ satisfies the required property.
\end{proof}

\section{Splitting necklaces}
\label{sec:necklace}

\subsection{The necklace theorem}
We turn now to the splitting necklace theorem. The first version of
this theorem, for two thieves, was proved by Goldberg and West
\cite{GoWe85}, and by Alon and West with a shorter proof
\cite{AlWe86}. The version with any number $q$ of thieves (but still
a multiple of $q$ of beads of each type) was proved by Alon
\cite{Al87}. The version proved here is slightly more general and
was proved by Alon, Moshkovitz and Safra \cite{AlMoSa06} (they
proved first a continuous version, and then proved that a
`rounding-procedure' is possible with flows).

 Suppose that the necklace has $n$ beads, each of
a certain type $i$, where $1\leq i\leq t$. Suppose there are $A_i$
beads of type $i$, $1\leq i\leq t$, $\sum_{i=1}^t A_i=n$, where
$A_i$ is not necessarily a multiple of $q$. A {\em $q$-splitting} of
the necklace is a partition of it into $q$ parts, each consisting of
a finite number of non-overlapping sub-necklaces of beads whose
union captures either $\lfloor A_i / q\rfloor$ or $\lceil A_i / q
\rceil$ beads of type $i$, for every $1\leq i\leq t$.

\begin{theorem}
Every necklace with $A_i$ beads of type $i$, $1\leq i\leq t$, has a
$q$-splitting requiring at most $t(q-1)$ cuts.
\end{theorem}

By a well-known trick (see \cite{Al87,Ma03}), it is enough to prove
it when $q=p$ is prime. The necklace is identified with the interval
$[0,n]$. The $k$th bead occupies uniformly the interval $[k-1,k]$.

\subsection{Encoding of the cuts -- $\mathsf{K}$ -- and of what get the thieves -- $\mathsf{L}$}
\label{subsec:encoding} We define a graph $\mathsf{S}:=(V,E)$ with
$1+pn$ vertices. This graph is seen as a 1-dimensional polytopal
complex. Define first the vertex set:
$$V:=\left\{(k,r):\,k\in\{1,\ldots,n\},r\in\{1,\ldots,p\}\right\}\cup\{o\}.$$
A vertex of $\mathsf{S}$ is either $o$ -- see it as the left origin
of the necklace -- or a position $k$ of a cut and a thief $r$.
Second, define the edges: $$E:=\left\{(k,r)(k',r'):\,|k-k'|=1\mbox{
and }r=r'\right\}\cup\left\{o(1,r):\,r\in\{1,\ldots,p\}\right\}.$$
$\mathsf{S}$ is then a collection of $p$ paths --
$P_1,P_2,\ldots,P_p$ -- having each $n$ edges, and a common vertex
$o$.

Consider the polytopal complex (actually, cubical complex)
$\mathsf{S}^{t(p-1)+1}$. Each vertex $\boldsymbol{v}$ is of the form
$(v_1,\ldots,v_{t(p-1)+1})$. Define $\mathsf{K}$ as the subcomplex
of $\mathsf{S}^{t(p-1)+1}$ such that always one of the $v_j$'s is of the form
$(n,r)$ with $r\in\{1,\ldots,p\}$. A vertex $\boldsymbol{v}$ of
$\mathsf{K}$ encodes a splitting of the necklace with at most
$t(p-1)$ cuts: each $v_j$ gives a cut and there are at most $t(p-1)$
such cuts that are different from $n$; now, consider the $j$th
sub-necklace and take the $v_{j'}=(k',r')$, with the smallest $j'$,
whose position $k'$ is on the right of this sub-necklace; the thief
$r'$ is the one who gets this sub-necklace.

\begin{figure}
\includegraphics[width=13cm]{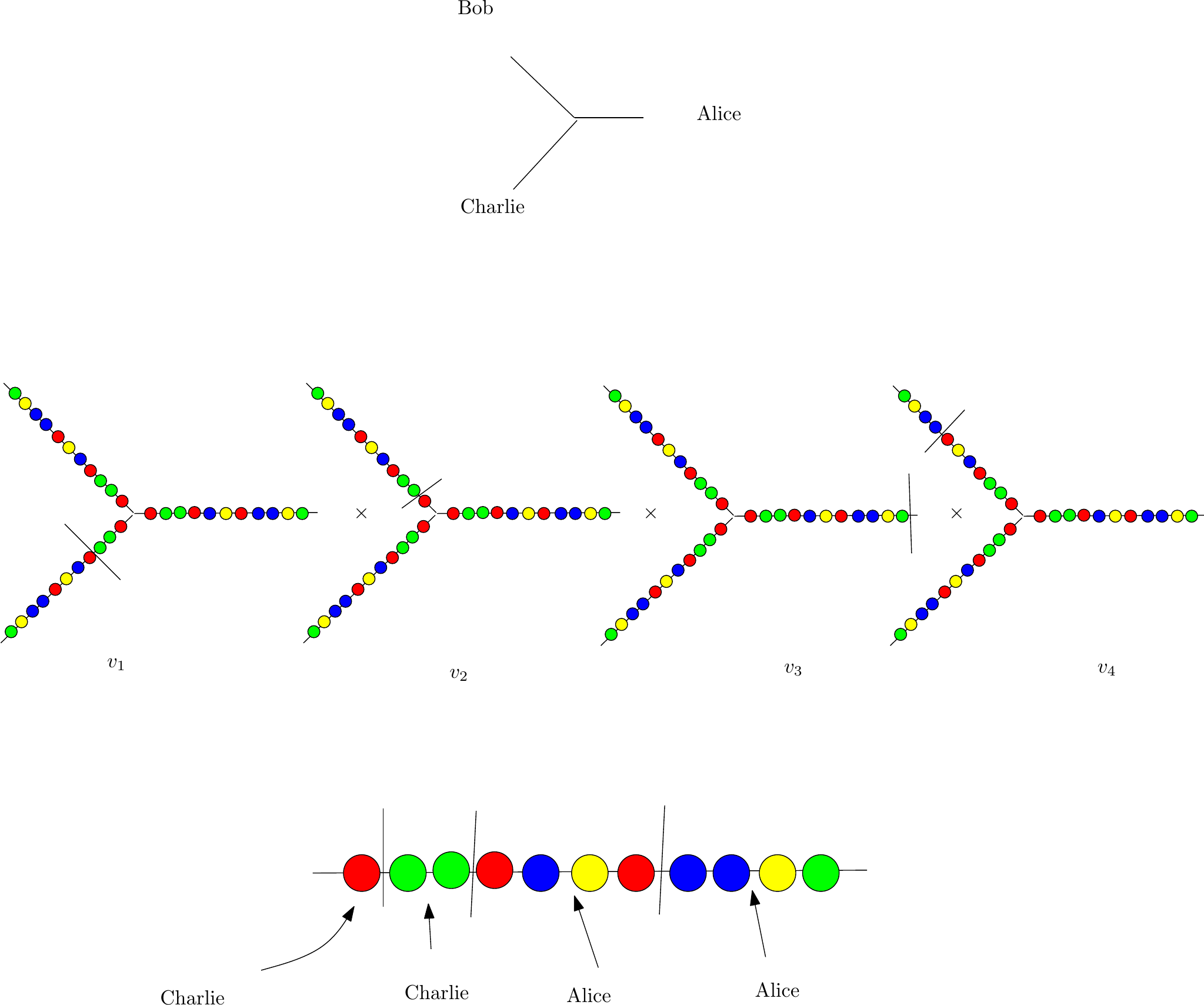}
\caption{Encoding of a splitting as a vertex of $\mathsf{K}$. Here
the $v_j$ of the form $(n,r)$ is $v_3$ with $r=\mbox{Bob}$.}
 \label{fig:codecuts}
\end{figure}

An example is given in Figure \ref{fig:codecuts}. Note that if the
first two cuts were exchanged, then the first subnecklace would get
to Bob and the second one still to Alice, the other assignments
remaining unchanged.

This encoding was proposed in the paper \cite{Me08} for the case
with two thieves.

Consider now the polytopal complex
$\mathsf{L}:=(\mathsf{\Delta}_{p-1})^t$, where
$\mathsf{\Delta}_{p-1}$ denotes the $(p-1)$-dimensional simplex,
whose vertices are $1,\ldots,p$. For each vertex $\boldsymbol{v}$ of
$\mathsf{K}$, define $\lambda_i(\boldsymbol{v})$ as the thief who
gets the largest amount of beads of type $i$ when one splits the
necklace according to $\boldsymbol{v}$ (using the positions of the
beads on the necklace as a total order, one can avoid a tie : in
case of equality, the thief with the lower position is considered as
advantaged). This thief is called the {\em $i$-winner}.

\subsection{Combinatorial polytopal Dold's theorem and proof of the necklace theorem}\label{subsec:combdold}
Consider now $\nu$ the cyclic shift $r\mapsto r+1$ modulo $p$. It
induces free actions on $\mathsf{K}$. The map
$$\begin{array}{cccc}\lambda: & \mathsf{K} & \rightarrow & \mathsf{L} \\ & \boldsymbol{v} & \mapsto &  \left(\lambda_1(\boldsymbol{v}),\ldots,\lambda_t(\boldsymbol{v})\right)
\end{array}$$
 is
then an equivariant polytopal map, that induces an equivariant chain
map $\mathcal{C}(\mathsf{K})\rightarrow\mathcal{C}(\mathsf{L})$ (one
uses here Theorem \ref{thm:chainmap} to derive the construction of this induced chain map).

\begin{lemma}\label{lem:polmap}
$\lambda$ is an equivariant polytopal map.
\end{lemma}

\begin{proof}
The equivariance is straightforward. Let us check that $\lambda$ is
polytopal, that is, that the image of a $d$-cube $\sigma$ in
$\mathsf{K}$ is included in a face of $\mathsf{L}$ of dimension at
most $d$.

Take a $d$-cube $\sigma$. It is defined by $d$ cuts, each of them
selected in one $\mathsf{S}$ in the product
$\mathsf{S}^{t(p-1)+1}$. Among these $d$ cuts, consider those sliding on type 1 beads on $\sigma$.  Denote by $d_1$ the number of such cuts. Recall that each cut is ``signed'', that is, indicates a thief $r$ (a cut is of the form $(k,r)$, where $k$ is a position and $r$ is a thief). Consider the hypergraph $H_1$ whose vertex set is
the set of thieves and whose hyperedges connect thieves that can get the same bead in the splitting encoded by the vertices of $\sigma$. Denote by $m_1$ the number of its hyperedges. We have then the following equality:
\begin{equation}\label{eq:in}\sum_{F\in H_1}|F|\leq d_1+m_1.\end{equation} Call $W$ the set of thieves that are $1$-winner among the vertices of $\sigma$.

Denote $a_1:=\left\lfloor\frac{A_1}{p}\right\rfloor$. Take the vertex of $\sigma$ that gives to the $1$-winner the
smallest number of beads of type $1$ among all the vertices of
$\sigma$. Let $w$ be this thief.

If $w$ gets $a_1$ beads, all the thieves get this number of beads.
Suppose then that $w$ gets $a_1+b$ beads where $b\geq 1$. Select in each edge of $H_1$ an {\em head}, which is precisely the thief who gets the corresponding bead in the present splitting. Let $n_x$ be the number of
thieves in $W$ who get $a_1+x$ beads of type 1. Now, consider a thief
$\neq w$ getting also exactly $a_1+b$ beads. There is at least one
hyperedge in which he is not the head, otherwise one can easily show, using the assumption
on $w$, that this thief cannot be a $1$-winner among the vertices of
$\sigma$.

Moreover, for a thief getting $a_1+x$ beads, there are at least
$b-x$ hyperedges in which he is not the head, since there exists a vertex of $\sigma$
giving him at least $a_1+b$ beads (he can be a $1$-winner).

Combining this two remarks, we get that
$n_b-1+\sum_{x<b}(b-x)n_x\leq\sum_{F\in H_1}(|F|-1)$, which implies according to equation (\ref{eq:in}) that
$n_b+\sum_{x<b}n_x\leq d_1+1$. Since $n_b+\sum_{x<b}n_x$ is the
number of vertices in $W$, i.e. the number of thieves that
may be $1$-winners, we get that the image of $\sigma$ by $\lambda$
in the first copy of $\mathsf{\Delta}_{p-1}$ in $\mathsf{L}$ is a
face of $\mathsf{\Delta}_{p-1}$ having at most $d_1+1$ vertices,
that is, being of dimension at most $d_1$.

The same holds for each type of beads. Finally, the dimension of the
minimal face containing $\lambda(\sigma)$ is at most
$d_1+d_2+\ldots+d_t=d$.
\end{proof}

\noindent{\bf Remark}: In the present case, $\lambda$ is a very particular polytopal map, since the image of (the vertex set of) any cube is the vertex set of a polytope (and not simply {\bf included} in the vertex set of polytope).

\begin{theorem}\label{thm:combdold}
Let $\lambda$ be an equivariant polytopal map $\mathsf{K}\rightarrow
\mathsf{L}$. Then $\mathsf{K}$ has a $t(p-1)$-face $\sigma$ whose
image by $\lambda$ is the $t(p-1)$-face of $\mathsf{L}$.
\end{theorem}

This theorem can be interpreted as a combinatorial polytopal Dold's
theorem. Dold's theorem is a theorem that generalizes Borsuk-Ulam
theorem for $\mathbb{Z}_p$ action, the Borsuk-Ulam theorem being the
case $p=2$ (see \cite{Do83,Ma03}).

Once this theorem is proved, we are done. Indeed, one has the
following lemma.

\begin{lemma}
If $\mathsf{K}$ has a $t(p-1)$-face $\sigma$ whose image by
$\lambda$ is the $t(p-1)$-face of $\mathsf{L}$, then $\sigma$ has at
least one vertex corresponding to a $p$-splitting.
\end{lemma}

\begin{proof}
Consider the hypergraphs $H_i$, and the numbers $d_i$ and $m_i$ as in the proof of Lemma
\ref{lem:polmap}. One of the $d_i$ is smaller than $p-1$ since
$\sum_{i=1}^td_i=t(p-1)$. Suppose w.l.o.g.
that it is $d_1$. Note that according to the assumption on
$\lambda$ and $\sigma$, each vertex of $H_1$ can be a $1$-winner for
a vertex of $\sigma$.

Defining similarly $n_x$ and $n_b$, as in the proof of Lemma
\ref{lem:polmap}. One has $n_b-1+\sum_{x<b}n_x(b-x)\leq p-1$. On the other hand,
there are $p$ vertices in $H_1$, thus $n_b+\sum_{x<b}n_x=p$.
Combining these two inequalities leads to
$n_b+\sum_{x<b}n_x(b-x)\leq n_b+\sum_{x<b}n_x$, that is
$\sum_{x<b}n_x(b-x)=\sum_{x<b}n_x$, which implies that $n_x=0$
whenever $x < b-1$. Thus, all thieves get $a_1$ or $a_1+1$ beads of
type 1. The beads of type 1 are fairly divided between the thieves.

Moreover, we see that $b=1$, and that $d_1$ is at
least $n_0+n_1-1=p-1$. It shows that $d_1=p-1$.
Hence, since we can make the same reasoning for all $d_i$
smaller than $p-1$, we get that all $d_i$ are equal to $p-1$,
and then that the division is fair for each type of beads.
\end{proof}

Hence, one has only to prove Theorem \ref{thm:combdold}. If we were
allowed to use homology, the proof will be a direct consequence of
the Hopf-Lefschetz formula. Moreover, such a prove would work by
contradiction. The following proof neither uses homology, nor works
by contradiction (it is constructive in a logical sense).

\begin{proof}[Proof of Theorem \ref{thm:combdold}]
Its proof uses three ingredients:
\begin{itemize}
\item[(i)] a sequence of $d$-chains $h_d$ in
$\mathcal{C}(\mathsf{K})$ such that $$h_0:=(o,(n,1),o,\ldots,o)\mbox{, }\,\,\partial
h_{2l+1}=\sum_{r=1}^p\nu^rh_{2l}\quad\mbox{and}\quad\partial
h_{2l+2}=(\nu-\nu^{-1})h_{2l+1}.$$
\item[(ii)] an equivariant chain map
$\eta_{\#}:\mathcal{C}(\mathsf{L})\rightarrow\mathcal{C}\left(\mathbb{Z}_p^{*(t(p-1))}\right)$,
where $*$ is the join operation.
\item[(iii)] a sequences of chain maps
$\phi_{d\#}:C_d\left(\mathbb{Z}_p^{*(t(p-1))}\right)\rightarrow\mathbb{Z}_p$
such that $\phi_{0\#}$ is equal to $1$ for a vertex in the first copy of $\mathbb{Z}_p$ in $\mathbb{Z}_p^{*(t(p-1))}$ and $0$ elsewhere,
$$\phi_{(2l+1)\#}\circ(\nu^{-1}-\nu)=\phi_{(2l)\#}\circ\partial\quad\mbox{and}\quad
\phi_{(2l+2)\#}\circ\left(\sum_{r=1}^p\nu^r\right)=\phi_{(2l+1)\#}\circ\partial.$$
\end{itemize}

\noindent (i) In each copy of $\mathsf{S}$, we have $p$ paths
$P_1,\ldots,P_p$ (defined above in Subsection
\ref{subsec:encoding}). We orient each $P_r$ them from $o$ to its
endpoint $(n,r)$. Define, with the notation introduced in Subsection
\ref{subsec:product}
\begin{eqnarray}
\tilde{h}_0 & := & o \\
\tilde{h}_{2l+1} & := & \sum_{r=1}^p\nu^r\tilde{h}_{2l}\otimes P_1 \\
\tilde{h}_{2l+2} & := & (\nu-\nu^{-1})\tilde{h}_{2l+1}\otimes P_1
\end{eqnarray}
Note that $\tilde{h}_d\in C_d(\mathsf{S}^d)$. Then define
$$h_d:=\partial(\tilde{h}_d\otimes P_1\otimes o\otimes\ldots\otimes o)+\tilde{h}_d\otimes o\otimes o\otimes\ldots\otimes o.$$
Using the fact that $(\nu-\nu^{-1})\circ(\sum_{r=0}^{p-1}\nu^r)=0$,
the checking that $(h_d)$ satisfies the required relation is
straightforward.

\medskip

\noindent (ii) Define $\mathsf{sd}(\mathsf{L})$ to be the
barycentric subdivision of $\mathsf{L}$. There is a natural chain
map $\mbox{sd}_{\#}$ that maps a face $(\sigma,\epsilon)$ of
$\mathsf{L}$ to the chain induced by $(\sigma,\epsilon)$ on
$\mathsf{sd}(\mathsf{L})$.

Then, consider the poset $\mathcal{P}$ of the faces of
$\mathsf{L}$, and let $\Delta(\mathcal{P})$ be the order complex of
$\mathcal{P}$, that is the simplicial complex whose vertices are the
elements of $\mathcal{P}$ and whose faces are the chains of
$\mathcal{P}$. One has the classical isomorphism
$\mathsf{sd}(\mathsf{L})\simeq\Delta(\mathcal{P})$.

Define a chain map $g_{\#}:
\mathcal{C}(\Delta(\mathcal{P}))\rightarrow\mathcal{C}\left(\mathbb{Z}_p^{*(t(p-1))}\right)$
by taking in each orbit of $\Delta(\mathcal{P})$ a face $\sigma$ of
$\mathsf{L}$ (recall that $\mathbb{Z}_p$ acts on $\mathcal{P}$) and
by defining $g(\sigma)$ to be any vertex in the
($\dim\sigma+1$)-copy of $\mathbb{Z}_p$ in
$\mathbb{Z}_p^{*(t(p-1))}$. The map $g$ is then an equivariant
simplicial map.

Finally define $\eta_{\#}:=g_{\#}\circ\mbox{sd}_{\#}$.

\medskip

\noindent (iii) This sequence can be constructed using methods presented in
\cite{Me07,HaSaScZi07}.

\bigskip

Then, we show by induction that
$$(\phi_{(2l)\#}\circ\eta_{\#}\circ\lambda_{\#})\left(\left(\sum_{r=1}^p\nu^r\right)h_{2l}\right)=(-1)^l\mbox{
mod
}p\quad\mbox{and}\quad(\phi_{(2l+1)\#}\circ\eta_{\#}\circ\lambda_{\#})\left((\nu-\nu^{-1})h_{2l+1}\right)=(-1)^{l+1}\mbox{
mod }p,$$ for $l=0,\ldots,\lfloor (t(p-1)-1)/2\rfloor$. Start with $l=0$.
One has
$(\phi_{0\#}\circ\eta_{\#}\circ\lambda_{\#})(\sum_{r=1}^p\nu^rh_0)=
\phi_{0\#}\left(\sum_{r=1}^p\nu^r\eta_{\#}(1,1,\ldots,1)\right)$, where $(1,1,\ldots,1)\in V((\mathsf{\Delta}_{p-1})^t)$. Hence, $(\phi_{0\#}\circ\eta_{\#}\circ\lambda_{\#})(\sum_{r=1}^p\nu^rh_0)=1$.
After that, the formulas above are proved by a straightforward induction.

Now, if $t(p-1)-1$ is even, we have
$(\phi_{(t(p-1)-1)\#}\circ\eta_{\#}\circ\lambda_{\#})\left(\left(\sum_{r=1}^p\nu^r\right)h_{t(p-1)-1}\right)\neq
0$. It can be rewritten
$(\phi_{(t(p-1)-1)\#}\circ\eta_{\#}\circ\lambda_{\#})(\partial
h_{t(p-1)})\neq 0$, or equivalently
$(\phi_{(t(p-1)-1)\#}\circ\partial\circ\eta_{\#})\lambda_{\#}(h_{t(p-1)})\neq
0$, which shows that $\lambda_{\#}(h_{t(p-1)})\neq 0$
 and hence that
there is an oriented $t(p-1)$-face $(\sigma,\epsilon)$ of
$\mathsf{K}$ whose image by $\lambda_{\#}$ is nonzero. The same holds if $t(p-1)-1$ is odd.
\end{proof}

\section{Discussion}

\label{sec:discussion}

\subsection{Cubical maps}

We have outlined the fact that our definition of a polytopal map,
when specialized for cubical complexes, does not lead to the cubical
map defined in its full generality (see Subsection
\ref{subsec:polymap}). Whether a more general version of a polytopal
map is possible, which would generalize correctly the notion of a
cubical map in its full generality is an open question.

An application might be a Sperner lemma for a polytopal complex that
generalizes simultaneously the classical Sperner's lemma and the
cubical Sperner's lemma found by Ky Fan \cite{Fa60}.

\begin{theorem}[(Classical) Sperner's lemma]
Let $\mathsf{K}$ be a simplicial complex that refines the
$d$-dimensional simplex $\mathsf{\Delta}_d$ and let
$\lambda:V(\mathsf{K})\rightarrow\mathsf{\Delta}_d$ be a labelling such
that $\lambda(v)$ is a vertex of the minimal face containing $v$.
Then there is a small $d$-simplex of $\mathsf{K}$ whose vertices are
one-to-one mapped on the vertices of $\mathsf{\Delta}_d$.
\end{theorem}

\begin{theorem}[Cubical Sperner's lemma]
Let $\mathsf{K}$ be a cubical subdivision of the $d$-dimensional
cube $\mathsf{\square}_d$ and let
$\lambda:V(\mathsf{K})\rightarrow\mathsf{\square}_d$ be a labelling
such that
\begin{enumerate}
\item $\lambda(v)$ is a vertex of the minimal face containing $v$,
\item two adjacent vertices are mapped by $\lambda$ on adjacent vertices or
the same vertex.
\end{enumerate}
Then there is a small $d$-cube of $\mathsf{K}$ whose vertices are
one-to-one mapped on the vertices of $\mathsf{\square}_d$.
\end{theorem}

\subsection{Approximation theorem}

A well-known fact about simplicial map is that if there is
continuous map from a simplicial complex $\mathsf{K}$ to another one
$\mathsf{L}$, then there is a subdivision $\mathsf{sd}^N\mathsf{K}$
and a simplicial map
$\lambda:\mathsf{sd}^N\mathsf{K}\rightarrow\mathsf{L}$ that is a
simplicial approximation of $f$.

The author is not aware of a similar property for cubical maps.

The {\em cubical subdivision} $\mathsf{csd}$ of a cube is the cartesian product of
the barycentric subdivision of the $1$-dimensional cubes whose
product defines the cube.

We make the following conjecture, which seems to be a first step
toward a cubical approximation theorem.

\begin{conjecture}
Let $P$ and $Q$ be two cubes, and let $f:V(P)\rightarrow V(Q)$. Then
there exists $N\geq 0$ and a cubical map $g:\mathsf{csd}^N
P\rightarrow Q$ such that for each $v\in V(P)$ one has $f(v)=g(v)$
(where $Q$ is identified with its face complex).
\end{conjecture}

As a starting point, let $P:=[0,1]$ and $Q:=[0,1]\times[0,1]$.
Define $f$ a continuous function mapping $0$ to $(0,0)$ and $1$ to
$(1,1)$. One has to subdivide $P$ ($N:=1$).

In the case of simplicial complex, such a situation never arises.

\subsection{Cup product}

When one defines the cohomology ring of a simplicial complex, one
needs to define the cup product of two cochains. Is it possible to
get the cup product of two cochains of a polytopal complex without
passing through a simplicial subdivision for instance?

\subsection{Algorithms}

Finally, the classical questions in this topic are:
\begin{itemize}
\item is there a constructive proof of the combinatorial polytopal Dold's theorem (Theorem \ref{thm:combdold})? Note that
from a purely logical point of view, our proof is constructive since
we do neither use the choice axiom nor use contradiction. But our
proof provides no algorithm.
\item is there a constructive proof of the necklace theorem?
\item is there a polynomial algorithm that solves the necklace
problem?
\end{itemize}

\bibliographystyle{amsplain}
\bibliography{Combinatorics}

\providecommand{\bysame}{\leavevmode\hbox to3em{\hrulefill}\thinspace}
\providecommand{\MR}{\relax\ifhmode\unskip\space\fi MR }
\providecommand{\MRhref}[2]{%
  \href{http://www.ams.org/mathscinet-getitem?mr=#1}{#2}
}
\providecommand{\href}[2]{#2}
\begin{thebibliography}{10}

\bibitem{Al87}
N.~Alon, \emph{Splitting necklaces}, Adv. in Math. \textbf{63} (1987),
  247--253.

\bibitem{AlMoSa06}
N.~Alon, D.~Moshkovitz, and S.~Safra, \emph{Algorithmic construction of sets
  for $k$-restrictions}, ACM Trans. Algorithms \textbf{2} (2006), 153--177.

\bibitem{AlWe86}
N.~Alon and D.~West, \emph{The {B}orsuk-{U}lam theorem and bisection of
  necklaces}, Proc. Amer. Math. Soc. \textbf{98} (1986), 623--628.

\bibitem{Do83}
A.~Dold, \emph{Simple proofs of some {B}orsuk-{U}lam results}, Contemp. Math.
  \textbf{19} (1983), 65--69.

\bibitem{EhHe95}
R.~Ehrenborg and G.~Hetyei, \emph{Generalizations of {B}axter's theorem and
  cubical homology}, J. Combin. Theory Ser. A \textbf{69} (1995), 233--287.

\bibitem{Fa60}
K.~Fan, \emph{Combinatorial properties of certain simplicial and cubical vertex
  maps}, Arch. Math. \textbf{11} (1960), 368--377.

\bibitem{GoWe85}
C.~H. Goldberg and D.~West, \emph{Bisection of circle colorings}, SIAM J.
  Algebraic Discrete Methods \textbf{6} (1985), 93--106.

\bibitem{HaSaScZi07}
B.~Hanke, R.~Sanyal, C.~Schultz, and G.~Ziegler, \emph{Combinatorial {S}tokes
  formulas via minimal resolutions}, J. of Combin. Theory, Series A (to
  appear).

\bibitem{Ma03}
J.~Matou\v{s}ek, \emph{Using the {B}orsuk-{U}lam theorem}, Springer Verlag,
  Berlin--Heidelberg--New York, 2003.

\bibitem{Me07}
F.~Meunier, \emph{A $\mathbb{Z}_q$-{F}an theorem}, Technical report, presented
  at the ``Topological combinatorics'' workshop, Stockholm (2006).

\bibitem{Me08}
\bysame, \emph{Discrete splittings of the necklace}, Mathematics of Operation
  Research (To appear).

\bibitem{Mu84}
J.~R. Munkres, \emph{Elements of algebraic topology}, Perseus Book Publishing,
  1984.

\bibitem{Zi02}
G.~Ziegler, \emph{Generalized {K}neser coloring theorems with combinatorial
  proofs}, Invent. Math. \textbf{147} (2002), 671--691.

\end{thebibliography}

\end{document}